\documentclass[oneside,english,11pt]{amsart}
\usepackage[T1]{fontenc}
\usepackage[latin9]{inputenc}
\usepackage{amstext}
\usepackage{amsthm}
\usepackage{amssymb}
\usepackage{graphicx}

\makeatletter

\providecommand{\tabularnewline}{\\}

\numberwithin{equation}{section}
\theoremstyle{plain}
\newtheorem{thm}{\protect\theoremname}
\theoremstyle{plain}
\newtheorem{prop}[thm]{\protect\propositionname}
\theoremstyle{plain}
\newtheorem{lem}[thm]{\protect\lemmaname}
\theoremstyle{definition}
\newtheorem{defn}[thm]{\protect\definitionname}

\makeatother

\usepackage{babel}
\providecommand{\definitionname}{Definition}
\providecommand{\lemmaname}{Lemma}
\providecommand{\propositionname}{Proposition}
\providecommand{\theoremname}{Theorem}

\newcommand\dx{\,\textrm{d}x}
\newcommand\dt{\,\textrm{d}t}
\newcommand\du{\,\textrm{d}u}
\newcommand\dv{\,\textrm{d}v}

\begin{document}
\title{Self-similar Differential Equations}
\author{Leon Q. Brin and Joe Fields}
\begin{abstract}
Differential equations where the graph of some derivative of a function
is composed of a finite number of similarity transformations of the
graph of the function itself are defined. We call these self-similar
differential equations (SSDEs) and prove existence and uniqueness
of solution under certain conditions. While SSDEs are not ordinary
differential equations, the technique for demonstrating existence
and uniqueness of SSDEs parallels that for ODEs. This paper appears
to be the first work on equations of this nature.
\end{abstract}

\maketitle

\section{Introduction}

\begin{figure}
\begin{centering}
(a) \raisebox{-.75in}{\includegraphics[width=1.75in]{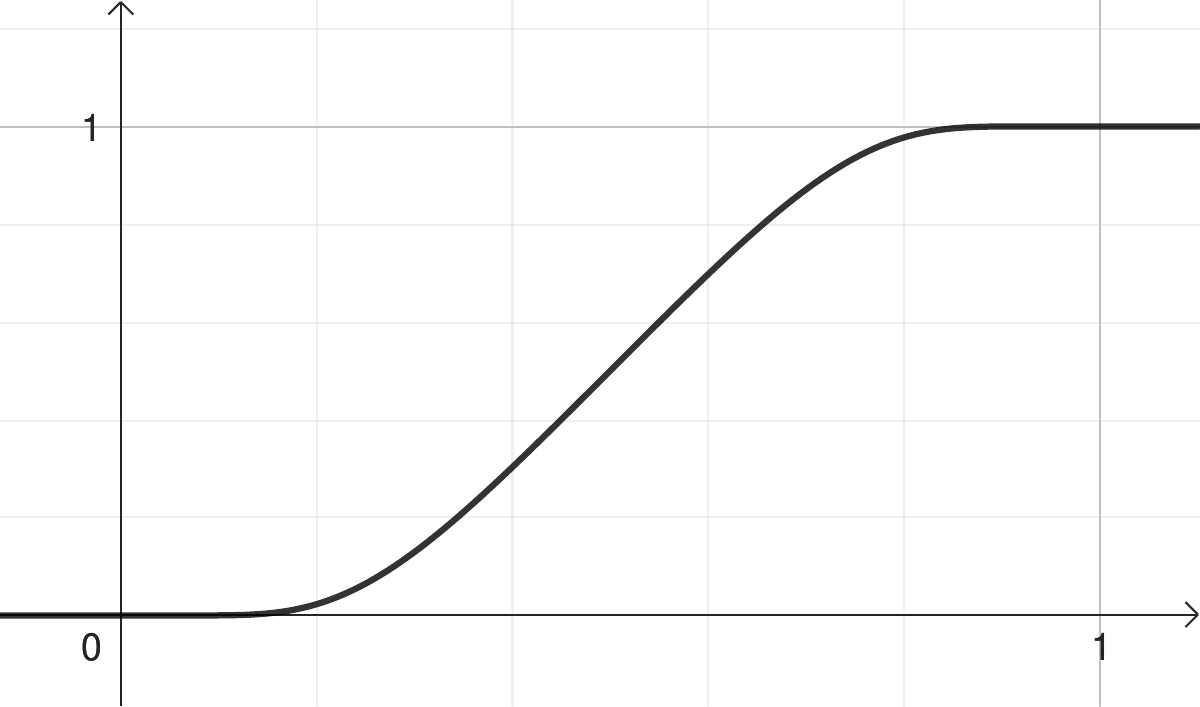}} 
\end{centering}

\vspace{.12in}

\begin{centering}
(b) \raisebox{-.75in}{\includegraphics[width=1.75in]{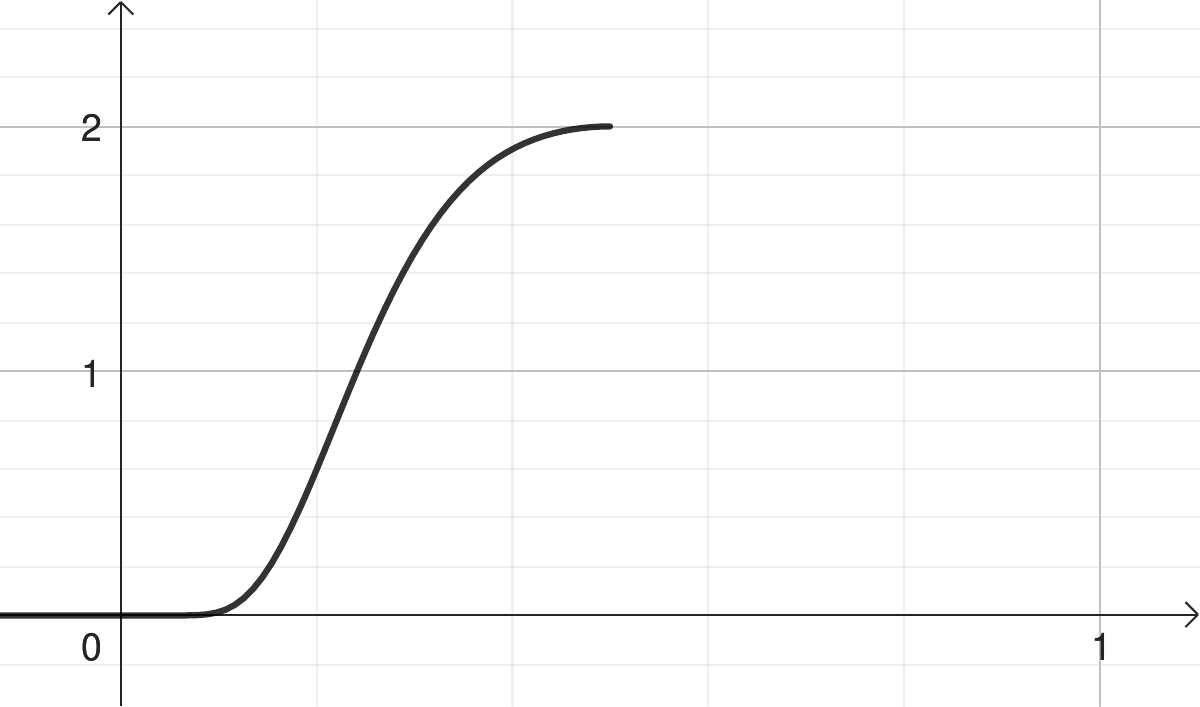}} 
\hspace{.5in}
(c) \raisebox{-.75in}{\includegraphics[width=1.75in]{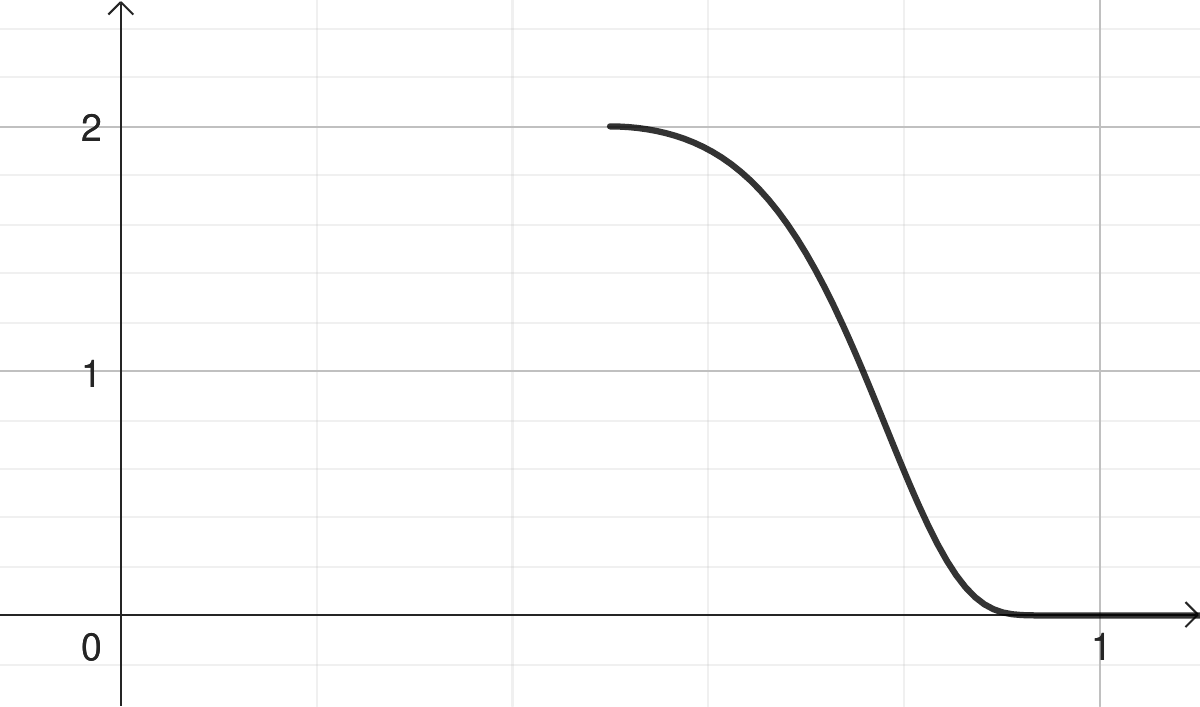}}
\end{centering}

\caption{(a) transition function $f$, (b) $f'$ on $(-\infty,1/2)$, (c) $f'$
on $(1/2,\infty)$ } \label{fig:transition}
\end{figure}

As motivation for upcoming definitions, consider the graph of the
smooth, monotonic transition function $y$ in Figure \ref{fig:transition}a,
which smoothly rises from the origin to $(1,1)$. Because it connects
two constant states, its derivative must be zero at both $0$ and
$1$. As drawn, it possesses a certain symmetry. The maximum value
of its derivative occurs at $1/2$ and $f'(1/2-\delta)=f'(1/2+\delta)$
for $\delta\in[0,1/2]$. Moreover, its derivative is increasing on
$(0,1/2)$ and decreasing on $(1/2,1)$. These observations suggest
that if such a function exists, its derivative over $[0,1/2]$ (Figure
\ref{fig:transition}b) looks like a smooth monotonic transition function
that rises smoothly from the origin to $(1/2,2)$, and its derivative
over $[1/2,1]$ (Figure \ref{fig:transition}c) looks like a smooth
monotonic transition function that falls smoothly from $(1/2,2)$
to $(1,0)$. In other words the graph of the derivative looks like
a patchwork of the graphs of two functions that each look a lot like
the whole function, giving it self-similarity in the first derivative.
This paper addresses the question of whether a function $f$ exists
where the graph of $f'$ over $[0,1/2]$ and the graph of $f'$ over
$[1/2,1]$ more than just look similar to $f$ in a vague sense, but
are similar in the mathematical sense. Imposing this idea on the function
we would necessarily have that the graph of $f'$ over $[0,1/2]$
be exactly the graph of $f$ stretched vertically by a factor of 2
and compressed horizontally by a factor of $1/2$. Likewise the graph
of $f'$ over $[1/2,1]$ would be exactly the graph of $f$ stretched
vertically by a factor of 2, compressed horizontally by a factor of
$1/2$, and reflected horizontally. To be more precise, it would necessarily
be that 
\begin{align}
f'(x) & =\begin{cases}
2f(2x) & 0\leq x\leq1/2\\
2f(2-2x) & 1/2<x\leq1
\end{cases}.\label{eq:transition}
\end{align}
Given that the function $f$ depicted in Figure \ref{fig:transition}
is defined by 
\[
f(x)=\frac{g(x)}{g(x)+g(1-x)}\quad\text{where}\quad g(x)=\begin{cases}
e^{-1/x} & x>0\\
0 & x\leq0
\end{cases}
\]
it is easy to verify that this particular $f$ does not satisfy \eqref{eq:transition}.
It will be shown in this paper that there does, however, exist a function
satisfying \eqref{eq:transition} and that functions with self-similarity
in their derivatives form a definable class of functions.

First order ordinary differential equations are given in the form
$f'(x)=g(f(x),x)$, and it is known, for example, that the initial
value problem $f'(x)=g(f(x),x)$; $f(x_{0})=f_{0}$ has a unique solution
provided, among other things, that $g$ is Lipschitz in $f$ \cite[pp.~106-113]{boyce}.
Equation \eqref{eq:transition} does not take the form $f'(x)=g(f(x),x)$
but rather $f'(x)=g(f(2x),f(2-2x))$, so the standard result does
not apply, and it is not even immediately clear what one would mean
by saying $g$ satisfied a Lipschitz condition in $f$. In short,
standard techniques of differential equations cannot be brought to
bear on \eqref{eq:transition} as the arguments on the left and right
hand sides differ. Nonetheless, we desire to prove existence and uniqueness
of the initial value problem defined by equation \eqref{eq:transition},
$f(0)=0$, on the interval $[0,1]$. By inspection, the trivial solution,
$f(x)=0$, is a solution, but there is no theory to suggest that this
solution is unique. Indeed it is not, and this will be proven presently.

Define

\begin{align*}
T(f)(x) & =\begin{cases}
2f(2x) & 0\leq x\leq1/2\\
2f(2-2x) & 1/2<x\leq1
\end{cases}\\
S(f)(x) & =\int_{0}^{x}T(f)(t)\dt\\
A & =\int_{0}^{1}f(x)\dx
\end{align*}
for any function $f$ integrable over $[0,1]$.

\begin{prop}
\label{prop:Sf-values}Given an integrable function $f$ on $[0,1]$,
$S(f)(1/2)=A$ and $S(f)(1)=2A$.
\end{prop}

\begin{proof}
Let $f$ be an integrable function defined on $[0,1]$. Then
\begin{align*}
S(f)\left(\frac{1}{2}\right) & =\int_{0}^{1/2}T(f)(t)\dt\\
 & =\int_{0}^{1/2}2f(2t)\dt=A
\end{align*}
and
\begin{align*}
S(f)\left(1\right) & =S(f)\left(\frac{1}{2}\right)+\int_{1/2}^{1}T(f)(t)\dt\\
 & =A+\int_{1/2}^{1}2f(2-2t)\dt=2A,
\end{align*}
both of which can be verified by simple change of variable.
\end{proof}

\begin{prop}
\label{prop:integralSf}Given an integrable function $f$ on $[0,1]$,
${\displaystyle \int_{0}^{1}S(f)(x)\dx=\int_{0}^{1}f(x)\dx}$.
\end{prop}

\begin{proof}
Let $f$ be an integrable function defined on $[0,1]$. Then
\begin{align*}
\int_{0}^{1}S(f)(x)\dx & =\int_{0}^{1/2}S(f)(x)\dx+\int_{1/2}^{1}S(f)(x)\dx\\
 & =\int_{0}^{1/2}\int_{0}^{x}T(f)(t)\dt\dx+\int_{1/2}^{1}\int_{0}^{x}T(f)(t)\dt\dx\\
 & =\int_{0}^{1/2}\int_{0}^{x}T(f)(t)\dt\dx+\int_{1/2}^{1}\left[\int_{0}^{1/2}T(f)(t)\dt+\int_{1/2}^{x}T(f)(t)\dt\right]\dx\\
 & =\int_{0}^{1/2}\int_{0}^{x}T(f)(t)\dt\dx+\int_{1/2}^{1}\int_{1/2}^{x}T(f)(t)\dt\dx+\int_{1/2}^{1}\int_{0}^{1/2}T(f)(t)\dt\dx\\
 & =\int_{0}^{1/2}\int_{0}^{x}2f(2t)\dt\dx+\int_{1/2}^{1}\int_{1/2}^{x}2f(2-2t)\dt\dx+\frac{1}{2}S(f)\left(\frac{1}{2}\right).
\end{align*}
By proposition 1, $\frac{1}{2}S(f)\left(\frac{1}{2}\right)=\frac{1}{2}A$,
and by direct calculation, $\int_{0}^{1/2}\int_{0}^{x}2f(2t)\dt\dx+\int_{1/2}^{1}\int_{1/2}^{x}2f(2-2t)\dt\dx=\frac{1}{2}A$,
completing the proof.
\end{proof}
Details of the calculation in proposition \ref{prop:integralSf} will
be demonstrated in generality later. In preparation for our existence
and uniqueness theorem, a lemma is needed.
\begin{lem}
\label{lem:max}If $f$ is continuous on $[L,R]$ and ${\displaystyle \int_{L}^{R}f(x)\dx=0}$
then 
\[
\max_{x\in[L,R]}\left|\int_{L}^{x}f(t)\dt\right|\leq\frac{1}{2}(R-L)\max_{x\in[L,R]}|f(x)|.
\]
\end{lem}

\begin{proof}
Suppose $f$ is continuous on $[L,R]$ and $\int_{L}^{R}f(x)\dx=0$.
Let $M=\max_{x\in[L,R]}|f(x)|$, and define $g(x)=\left|\int_{L}^{x}f(t)\dt\right|$
for $x\in[L,R]$. Because $g$ is continuous on $[L,R]$ and $g(L)=g(R)=0$
there must exist $c\in(L,R)$ such that $g(c)=\max_{x\in[L,R]}g(x)$.
Hence
\[
\max_{x\in[L,R]}g(x)=g(c)=\left|\int_{L}^{c}f(t)\dt\right|\leq\int_{L}^{c}\left|f(t)\right|\dt\leq(c-L)\max_{x\in[L,c]}|f(x)|\leq(c-L)M.
\]
Because $\int_{L}^{R}f(x)\dx=0$, $\int_{c}^{R}f(t)\dt=-\int_{L}^{c}f(t)\dt$,
from which it follows
\[
\max_{x\in[L,R]}g(x)=\left|\int_{c}^{R}f(t)\dt\right|\leq\int_{c}^{R}\left|f(t)\right|\dt\leq(R-c)\max_{x\in[c,R]}|f(x)|\leq(R-c)M.
\]
Because $c$ lies between $L$ and $R$, $c-L\leq\frac{1}{2}(R-L)$
or $R-c\leq\frac{1}{2}(R-L)$. Either way, this completes the proof.
\end{proof}
Application of lemma \ref{lem:max} and the contraction mapping principle,
alternatively unnamed or called the contraction mapping theorem ,
\cite[pp.~283-284]{depree} \cite[p.~137]{browder} \cite[p.~98]{carothers}
will provide an existence and uniqueness result for \eqref{eq:transition}.
\begin{thm}
\label{thm:existence_uniqeness_1}For each real value, $A$, there
exists a unique solution, $f$, over $[0,1]$ of \eqref{eq:transition}
with $f(0)=0$ and ${\displaystyle \int_{0}^{1}f(x)\dx=A}$.
\end{thm}

\begin{figure}
\begin{centering}
(a) \raisebox{-.75in}{\includegraphics[width=1.75in]{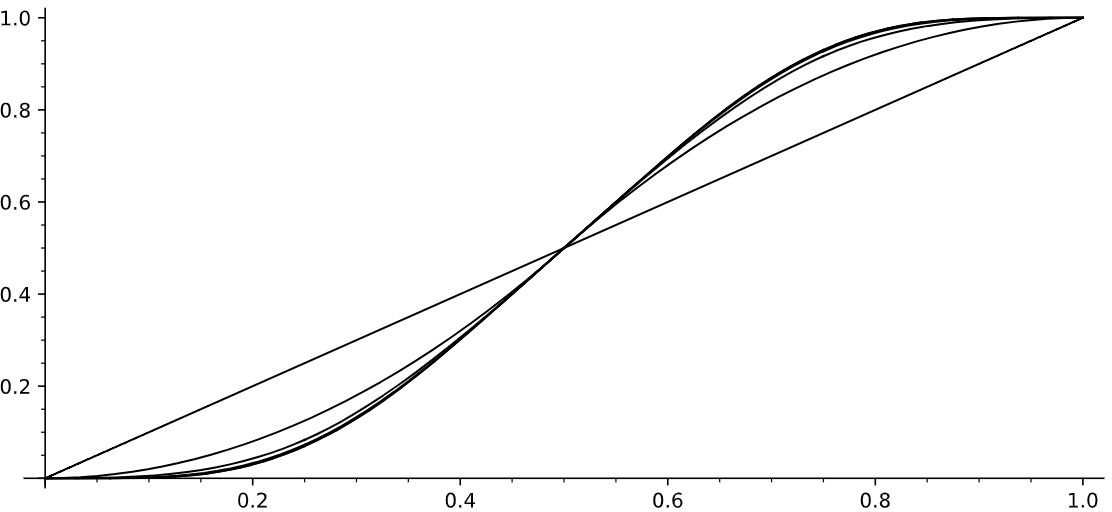}}
(b) \raisebox{-.75in}{\includegraphics[width=1.75in]{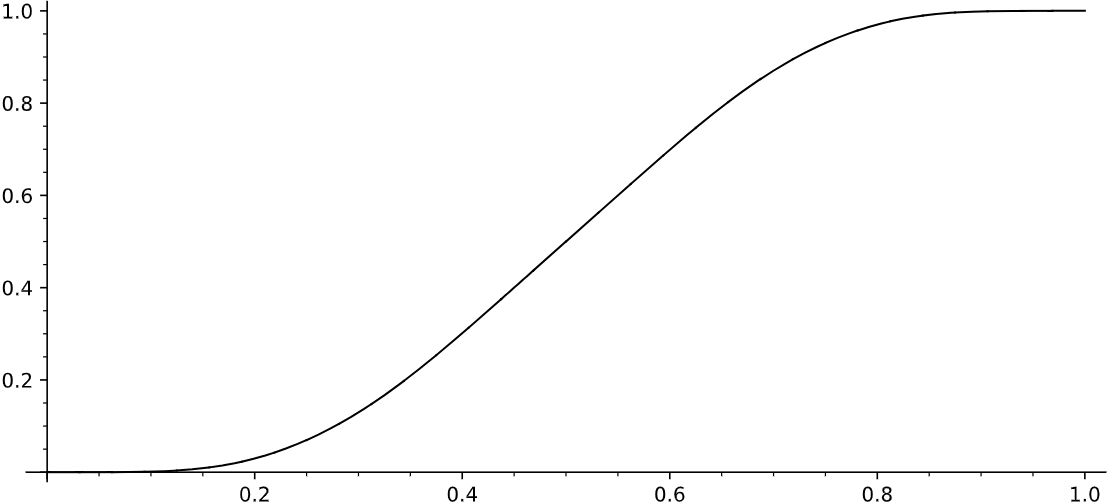}
}
\par\end{centering}
\caption{(a) iterates of $S$ beginning with $f_{0}(x)=x$, (b) the fifth iteration} \label{fig:transition_graph}
\end{figure}

\begin{proof}
Let $A$ be any real number and let $g$ and $h$ be continuous functions
such that ${\displaystyle \int_{0}^{1}g(x)\dx=\int_{0}^{1}h(x)\dx=A}$.
We wish to compare ${\displaystyle \max_{x\in[0,1]}\left|S(g)(x)-S(h)(x)\right|}$
with ${\displaystyle \max_{x\in[0,1]}\left|g(x)-h(x)\right|}$. To
that end,
\begin{align*}
S(g)(x)-S(h)(x) & =\int_{0}^{x}T(g)(t)\dt-\int_{0}^{x}T(h)(t)\dt=\int_{0}^{x}\left(T(g)(t)-T(h)(t)\right)\dt\\
 & =\begin{cases}
\int_{0}^{x}\left(2g(2t)-2h(2t)\right)\dt & 0\leq x\leq1/2\\
S(g)\left(\frac{1}{2}\right)-S(h)\left(\frac{1}{2}\right)+\int_{1/2}^{x}\left(2g(2-2t)-2h(2-2t)\right)\dt & 1/2<x\leq1
\end{cases}
\end{align*}
By proposition \ref{prop:Sf-values}, $S(g)\left(\frac{1}{2}\right)=S(h)\left(\frac{1}{2}\right)$.
It follows that
\[
S(g)(x)-S(h)(x)=\begin{cases}
\int_{0}^{2x}\left(g(u)-h(u)\right)\du & 0\leq x\leq1/2\\
\int_{1}^{2-2x}\left(g(u)-h(u)\right)\du & 1/2<x\leq1
\end{cases}
\]
and therefore
\[
\max_{x\in[0,1/2]}\left|S(g)(x)-S(h)(x)\right|=\max_{x\in[0,1/2]}\left|\int_{0}^{2x}\left(g(u)-h(u)\right)\du\right|=\max_{x\in[0,1]}\left|\int_{0}^{x}\left(g(u)-h(u)\right)\du\right|
\]
and
\[
\max_{x\in[1/2,1]}\left|S(g)(x)-S(h)(x)\right|=\max_{x\in[1/2,1]}\left|\int_{1}^{2-2x}\left(g(u)-h(u)\right)\du\right|=\max_{x\in[0,1]}\left|\int_{0}^{x}\left(g(u)-h(u)\right)\du\right|.
\]
We conclude that
\[
\max_{x\in[0,1]}\left|S(g)(x)-S(h)(x)\right|=\max_{x\in[0,1]}\left|\int_{0}^{x}\left(g(u)-h(u)\right)\du\right|.
\]
But ${\displaystyle \int_{0}^{1}\left(g(u)-h(u)\right)\du=0}$, so
by lemma \ref{lem:max}
\[
\max_{x\in[0,1]}\left|\int_{0}^{x}\left(g(u)-h(u)\right)\du\right|\leq\frac{1}{2}\max_{x\in[0,1]}\left|g(x)-h(x)\right|.
\]

By proposition \ref{prop:integralSf}, $Z_{A}=\{f:[0,1]\to\mathbb{R}|f\text{ is continuous and }\int_{0}^{1}f(x)\dx=A\}$
is closed under $S$. By the above calculation, $S$ is contractive
on $Z_{A}$ under the max norm. Since $Z_{A}$ is complete with respect
to uniform convergence, the contraction mapping principle guarantees
a unique fixed point of $S$ in $Z$. This fixed point, then, is the
unique solution of the equation $f=S(f)$, the integral equation equivalent
of \eqref{eq:transition} with $f(0)=0$, within $Z_{A}$.
\end{proof}
Beyond the guarantee of uniqueness, the contraction mapping principle
provides an algorithm for approximating the solution $f$ for any
value $A$. For example, set $f_{0}(x)=x$, which gives $A=\int_{0}^{1}f_{0}(x)\dx=1/2$.
Then set $f_{1}(x)=S(f_{0})(x)$, $f_{2}(x)=S(f_{1})(x)$, and so
on to produce a sequence of functions that, in the limit, yield a
transition function $f$ satisfying \eqref{eq:transition} with $f(0)=0$
and $f(1)=1$. The results of the first five iterations are shown
in Figure \ref{fig:transition_graph}. The analogous procedure in
the theory of ordinary differential equations is most often referred
to as Picard iteration \cite[p.~106]{boyce}.

 Despite the differences between \eqref{eq:transition} and ordinary
differential equations, it is the same principle that provides proof
of existence and uniqueness of solutions.

Transition functions such as this are useful for defining bump functions
\cite[p.~41]{hirsch}. If we extend $f'$ to be zero outside the interval
$[0,1]$, and call this extension $\hat{f'}$, then $\hat{f'}$ is
a bump function with support $[0,1]$ and the product $\hat{f'}\left(\frac{x-a}{b-a}\right)\hat{f'}\left(\frac{d-x}{d-c}\right)$
is a bump function with support $[a,d]$ for any real numbers $a<b<c<d$.

\section{Generalization}

This section is concerned with defining a general class of differential
equations of which \eqref{eq:transition} is but one example and generalizing
the results presented in the introduction. Consistent with the motivational
example, definitions are taken from a geometric viewpoint. As such
we make no distinction between a function and its graph.

\subsection{Definitions}
\begin{defn}
We define the following terms.
\end{defn}

\begin{enumerate}
\item For $D_{1},D_{2}\subseteq\mathbb{R}$, the transformation $T:D_{1}\times\mathbb{R}\to D_{2}\times\mathbb{R}$
is called \textbf{function preserving} if whenever $g$ is a function
on $D_{1}$, $T(g)$ is a function on $D_{2}$.
\item If $x_{0}<x_{1}<\cdots<x_{n}$ and $T_{1},T_{2},\ldots,T_{n}$ are
function preserving maps $T_{j}:[x_{0},x_{n}]\times\mathbb{R}\to[x_{j-1},x_{j}]\times\mathbb{R}$,
then $\{T_{j}|j=1,2,\ldots,n\}$ is a \textbf{piecemealing} on $[x_{0},x_{n}]$.
\item Suppose $P=\{T_{j}\}$ is a piecemealing on $I=[x_{0},x_{n}]$ and
$y:I\to\mathbb{R}$. We define
\begin{enumerate}
\item $P(y):I\to\mathbb{R}$ by union: $P(y)=\cup_{j}T_{j}(y)$ and
\item $\mathcal{F}_{C}(y):I\to\mathbb{R}$ by $\mathcal{F}_{C}(y)(x)=C+\int_{x_{0}}^{x}P(y)(t)\dt$.
\end{enumerate}
\end{enumerate}
We now consider differential equations of the form
\begin{align}
y' & =P(y)\label{eq:first_order}\\
y(x_{0}) & =y_{0}\nonumber 
\end{align}
where $P$ is a piecemealing on $I=[x_{0},x_{n}]$. We refer to \eqref{eq:first_order}
as a \textbf{self-similar differential equation} or SSDE. Note that
equation \eqref{eq:transition} with $f(0)=0$ is equivalent to the
SSDE \eqref{eq:first_order} with $x_{0}=y_{0}=0$ and 
\[
P=\left\{ T_{1}\left[\begin{array}{c}
x\\
y
\end{array}\right]=\left[\begin{array}{cc}
\frac{1}{2} & 0\\
0 & 2
\end{array}\right]\left[\begin{array}{c}
x\\
y
\end{array}\right]+\left[\begin{array}{c}
0\\
0
\end{array}\right],\ T_{2}\left[\begin{array}{c}
x\\
y
\end{array}\right]=\left[\begin{array}{cc}
-\frac{1}{2} & 0\\
0 & 2
\end{array}\right]\left[\begin{array}{c}
x\\
y
\end{array}\right]+\left[\begin{array}{c}
1\\
0
\end{array}\right]\right\} .
\]
It is a simple matter to verify that transformations of the form 
\[
\left[\begin{array}{cc}
a & 0\\
c & d
\end{array}\right]\left[\begin{array}{c}
x\\
y
\end{array}\right]+\left[\begin{array}{c}
e\\
f
\end{array}\right]
\]
with $a\neq0$ are function preserving. 

\subsection{General results}

Given $x_{0}<x_{1}<\cdots<x_{n}$, we now consider SSDE's with piecemealings
of the form 
\begin{equation}
P=\left\{ T_{i}\left[\begin{array}{c}
x\\
y
\end{array}\right]=\left[\begin{array}{cc}
a_{i} & 0\\
0 & d_{i}
\end{array}\right]\left[\begin{array}{c}
x\\
y
\end{array}\right]+\left[\begin{array}{c}
e_{i}\\
f_{i}
\end{array}\right]|i=1,2,\ldots,n\right\} \label{eq:linear_piecemealing_FDE}
\end{equation}
 on $[x_{0},x_{n}]$ where $T_{i}([x_{0},x_{n}]\times\mathbb{R})=[x_{i-1},x_{i}]\times\mathbb{R}$
for each $i$.
\begin{prop}
\label{prop:basics1}Let $P$ be a piecemealing of the form \eqref{eq:linear_piecemealing_FDE}.
Then
\[
(a_{i},e_{i})=\begin{cases}
\left(\frac{x_{i}-x_{i-1}}{x_{n}-x_{0}},x_{i}-a_{i}x_{n}\right) & \text{for }a_{i}>0\\
\left(\frac{x_{i-1}-x_{i}}{x_{n}-x_{0}},x_{i}-a_{i}x_{0}\right) & \text{for }a_{i}<0
\end{cases}
\]
\end{prop}

\begin{proof}
Because $T_{i}$ maps vertical lines to vertical lines, it must be
that the image of the line $x=x_{0}$ is the line $x=x_{i-1}$ for
a temporal preserving transformation ($a_{i}>0$) and must be the
line $x=x_{i}$ for a temporal reversing transformation ($a_{i}<0$).
Similarly the image of the line $x=x_{n}$ must be $x=x_{i}$ for
a temporal preserving transformation and must be $x=x_{i-1}$ for
a temporal reversing transformation. The result follows by direct
calculation.
\end{proof}
\begin{prop}
\label{prop:basics2}Let $P$ be a piecemealing of the form \eqref{eq:linear_piecemealing_FDE}
and let $y:[x_{0},x_{n}]\to\mathbb{R}$ be integrable. Then for each
$k=0,1,2,\ldots,n$,
\[
\int_{x_{0}}^{x_{k}}P(y)(x) \dx=\sum_{i=1}^{k}\left[|a_{i}|d_{i}\int_{x_{0}}^{x_{n}}y(x)\dx+f_{i}(x_{i}-x_{i-1})\right].
\]
\end{prop}

\begin{proof}
Let $C=\sum_{i=1}^{k}f_{i}(x_{i}-x_{i-1})$, and let $u=\frac{x-e_{i}}{a_{i}}$.
Then
\begin{align*}
\int_{x_{0}}^{x_{k}}P(y)(x)\dx & =\sum_{i=1}^{k}\int_{x_{i-1}}^{x_{i}}P(y)(x)\dx\\
 & =\sum_{i=1}^{k}\int_{x_{i-1}}^{x_{i}}T_{i}(y)(x)\dx\\
 & =\sum_{i=1}^{k}\int_{x_{i-1}}^{x_{i}}\left[d_{i}y\left(\frac{x-e_{i}}{a_{i}}\right)+f_{i}\right]\dx\\
 & =C+\sum_{i=1}^{k}a_{i}d_{i}\begin{cases}
\int_{x0}^{x_{n}}y(u)\du & \text{if }a_{i}>0\\
\int_{x_{n}}^{x_{0}}y(u)\du & \text{if }a_{i}<0
\end{cases}\\
 & =C+\sum_{i=1}^{k}|a_{i}|d_{i}\int_{x_{0}}^{x_{n}}y(u)\du
\end{align*}
\end{proof}
\begin{prop}
\label{prop:basics3}Let $P$ be a piecemealing of the form \eqref{eq:linear_piecemealing_FDE}.
Let $y:[x_{0},x_{n}]\to\mathbb{R}$ be integrable, and let $Y(x)$
be the antiderivative of $y$ with $Y(x_{0})=0$. Then

\[
\int_{x_{0}}^{x_{n}}\int_{x_{0}}^{x}P(y)(t)\dt\dx=K+L-M+N
\]
where 
\[
\begin{gathered}K=\sum_{i=1}^{n}\int_{x_{i-1}}^{x_{i}}\int_{x_{0}}^{x_{i-1}}P(y)(t)\dt\dx\qquad L=\int_{x_{0}}^{x_{n}}Y(v)\dv\sum_{i=1}^{n}\frac{a_{i}^{3}d_{i}}{|a_{i}|}\\
M=(x_{n}-x_{0})\int_{x_{0}}^{x_{n}}y(x)\dx\sum_{i\ni a_{i}<0}a_{i}^{2}d_{i}\qquad N=\frac{1}{2}(x_{n}-x_{0})^{2}\sum_{i=1}^{n}a_{i}^{2}f_{i}
\end{gathered}
.
\]
\end{prop}

\begin{proof}
First,
\begin{align*}
\int_{x_{0}}^{x_{n}}\int_{x_{0}}^{x}P(y)(t)\dt\dx & =\sum_{i=1}^{n}\int_{x_{i-1}}^{x_{i}}\int_{x_{0}}^{x}P(y)(t)\dt\dx\\
 & =\sum_{i=1}^{n}\int_{x_{i-1}}^{x_{i}}\left[\int_{x_{0}}^{x_{i-1}}P(y)(t)\dt+\int_{x_{i-1}}^{x}P(y)(t)\dt\right]\dx\\
 & =K+\sum_{i=1}^{n}\int_{x_{i-1}}^{x_{i}}\int_{x_{i-1}}^{x}P(y)(t)\dt\dx
\end{align*}
But
\begin{align*}
\sum_{i=1}^{n}\int_{x_{i-1}}^{x_{i}}\int_{x_{i-1}}^{x}P(y)(t)\dt\dx & =\sum_{i=1}^{n}\int_{x_{i-1}}^{x_{i}}\int_{x_{i-1}}^{x}T_{i}(y)(t)\dt\dx\\
 & =\sum_{i=1}^{n}\int_{x_{i-1}}^{x_{i}}\int_{x_{i-1}}^{x}\left[d_{i}y\left(\frac{t-e_{i}}{a_{i}}\right)+f_{i}\right]\dt\dx
\end{align*}
Now let $w=x_{n}-x_{0}$ and note that $\sum_{i=1}^{n}\int_{x_{i-1}}^{x_{i}}\int_{x_{i-1}}^{x}f_{i}\dt\dx=\frac{1}{2}\sum_{i=1}^{n}f_{i}(x_{i}-x_{i-1})^{2}=\frac{1}{2}\sum_{i=1}^{n}f_{i}\left(a_{i}w\right)^{2}=N$.
It remains to show $\sum_{i=1}^{n}\int_{x_{i-1}}^{x_{i}}\int_{x_{i-1}}^{x}d_{i}y\left(\frac{t-e_{i}}{a_{i}}\right)\dt\dx=L-M$.
Letting $u=\frac{t-e_{i}}{a_{i}}$ and $v=\frac{x-e}{a_{i}}$
\begin{align*}
\sum_{i=1}^{n}\int_{x_{i-1}}^{x_{i}}\int_{x_{i-1}}^{x}d_{i}y\left(\frac{t-e_{i}}{a_{i}}\right)\dt\dx & =\sum_{i=1}^{n}a_{i}^{2}d_{i}\int_{\frac{x_{i-1}-e_{i}}{a_{i}}}^{\frac{x_{i}-e_{i}}{a_{i}}}\int_{\frac{x_{i-1}-e_{i}}{a_{i}}}^{v}y\left(u\right)\du\dv\\
 & =\sum_{i=1}^{n}\begin{cases}
a_{i}^{2}d_{i}\int_{x_{0}}^{x_{n}}\int_{x_{0}}^{v}y\left(u\right)\du\dv & \text{if }a_{i}>0\\
a_{i}^{2}d_{i}\int_{x_{n}}^{x_{0}}\int_{x_{n}}^{v}y\left(u\right)\du\dv & \text{if }a_{i}<0
\end{cases}
\end{align*}
But $\int_{x_{0}}^{x_{n}}\int_{x_{0}}^{v}y\left(u\right)\du\dv=\int_{x_{0}}^{x_{n}}\left(Y(v)-Y(x_{0})\right)\dv=\int_{x_{0}}^{x_{n}}Y(v)\dv$
since $Y(x_{0})=0$. Moreover $\int_{x_{n}}^{x_{0}}\int_{x_{n}}^{v}y\left(u\right)\du\dv=\int_{x_{n}}^{x_{0}}\left(Y(v)-Y(x_{n})\right)\dv=\int_{x_{n}}^{x_{0}}Y(v)\dv-wY(x_{n})$.
Hence 
\begin{align*}
\sum_{i=1}^{n}\int_{x_{i-1}}^{x_{i}}\int_{x_{i-1}}^{x}d_{i}y\left(\frac{t-e_{i}}{a_{i}}\right)\dt\dx & =\sum_{i=1}^{n}\begin{cases}
a_{i}^{2}d_{i}\int_{x_{0}}^{x_{n}}Y(v)\dv & \text{if }a_{i}>0\\
-a_{i}^{2}d_{i}\left[\int_{x_{0}}^{x_{n}}Y(v)\dv+wY(x_{n})\right] & \text{if }a_{i}<0
\end{cases}\\
 & =\sum_{i=1}^{n}\frac{a_{i}^{3}d_{i}}{|a_{i}|}\int_{x_{0}}^{x_{n}}Y(v)\dv-w\int_{x_{0}}^{x_{n}}y(x)\dx\sum_{i\ni a_{i}<0}a_{i}^{2}d_{i}\\
 & =L-M.
\end{align*}
\end{proof}
Inspired by propositions \ref{prop:basics2} and \ref{prop:basics3},
we will make extensive use of the linear system of $n+1$ equations
\begin{align}
y_{k} & =y_{0}+\sum_{i=1}^{k}\left[|a_{i}|d_{i}A+f_{i}(x_{i}-x_{i-1})\right]\quad k=1,2,\ldots,n\label{eq:linear_system_part1}\\
A & =y_{0}(x_{n}-x_{0})+\sum_{i=1}^{n}(y_{i-1}-y_{0})(x_{i}-x_{i-1})\label{eq:linear_system_part2}\\
 & \phantom{=y_{0}(x_{n}-x_{0})}-(x_{n}-x_{0})A\sum_{i\ni a_{i}<0}a_{i}^{2}d_{i}+\frac{1}{2}(x_{n}-x_{0})^{2}\sum_{i=1}^{n}a_{i}^{2}f_{i}\nonumber 
\end{align}
where it is understood that $a_{i},d_{i},f_{i},x_{i},y_{0}$ are known
quantities and $y_{1},y_{2},\ldots,y_{n},A$ are unknown quantities.
\begin{lem}
\label{lem:vals_of_y}Let $P$ be a piecemealing of the form \eqref{eq:linear_piecemealing_FDE}
such that $\sum_{i=1}^{n}\frac{a_{i}^{3}d_{i}}{|a_{i}|}=0$. If $y$
is a solution of \eqref{eq:first_order} on $[x_{0},x_{n}]$, then
\[
y(x_{k})=y_{0}+\sum_{i=1}^{k}\left[|a_{i}|d_{i}\int_{x_{0}}^{x_{n}}y(x)\dx+f_{i}(x_{i}-x_{i-1})\right].
\]
\end{lem}

\begin{proof}
Let $y$ be a solution of \eqref{eq:first_order} on $[x_{0},x_{n}]$.
By proposition \eqref{prop:basics2}, 
\[
\int_{x_{0}}^{x_{k}}P(y)(x)\dx=\sum_{i=1}^{k}\left[|a_{i}|d_{i}\int_{x_{0}}^{x_{n}}y(x)\dx+f_{i}(x_{i}-x_{i-1})\right].
\]
 Since $y$ is a solution of \eqref{eq:first_order}, it is also true
that
\[
\int_{x_{0}}^{x_{k}}P(y)(x)\dx=\int_{x_{0}}^{x_{k}}y'(x)\dx=y(x_{k})-y_{0},
\]
so we have
\[
y(x_{k})=y_{0}+\sum_{i=1}^{k}\left[|a_{i}|d_{i}\int_{x_{0}}^{x_{n}}y(x)\dx+f_{i}(x_{i}-x_{i-1})\right].
\]
\end{proof}
\begin{thm}
\label{thm:FDE_soln_is_lin_sys_soln}Let $P$ be a piecemealing of
the form \eqref{eq:linear_piecemealing_FDE} such that $\sum_{i=1}^{n}\frac{a_{i}^{3}d_{i}}{|a_{i}|}=0$.
If $y$ is a solution of \eqref{eq:first_order} on $[x_{0},x_{n}]$,
then 
\begin{align}
y_{k} & =y(x_{k})\quad k=1,2,\ldots,n\nonumber \\
A & =\int_{x_{0}}^{x_{n}}y(x)\dx\label{eq:linear_system_solution}
\end{align}
is a solution of system \eqref{eq:linear_system_part1}, \eqref{eq:linear_system_part2}.
\end{thm}

\begin{proof}
Let $y$ be a solution of \eqref{eq:first_order} on $[x_{0},x_{n}]$.
By lemma \ref{lem:vals_of_y}, 
\[
y(x_{k})=y_{0}+\sum_{i=1}^{k}\left[|a_{i}|d_{i}\int_{x_{0}}^{x_{n}}y(x)\dx+f_{i}(x_{i}-x_{i-1})\right],
\]
demonstrating that equations \eqref{eq:linear_system_part1} are satisfied
by \eqref{eq:linear_system_solution}. By proposition \ref{prop:basics3}
and the fact that $\sum_{i=1}^{n}\frac{a_{i}^{3}d_{i}}{|a_{i}|}=0$,
\begin{align*}
\int_{x_{0}}^{x_{n}}\int_{x_{0}}^{x}P(y)(t)\dt\dx & =\sum_{i=1}^{n}\int_{x_{i-1}}^{x_{i}}\int_{x_{0}}^{x_{i-1}}P(y)(t)\dt\dx\\
 & \phantom{=}-(x_{n}-x_{0})\int_{x_{0}}^{x_{n}}y(x)\dx\sum_{i\ni a_{i}<0}a_{i}^{2}d_{i}+\frac{1}{2}(x_{n}-x_{0})^{2}\sum_{i=1}^{n}a_{i}^{2}f_{i}.
\end{align*}
Again using the fact that $y$ is a solution of \eqref{eq:first_order},
\begin{align*}
\int_{x_{0}}^{x_{n}}\int_{x_{0}}^{x}P(y)(t)\dt\dx & =\int_{x_{0}}^{x_{n}}\int_{x_{0}}^{x}y'(t)\dt\dx=\int_{x_{0}}^{x_{n}}\left(y(x)-y(x_{0})\right)\dx\\
 & =\int_{x_{0}}^{x_{n}}y(x)\dx-y_{0}(x_{n}-x_{0})
\end{align*}
and 
\begin{align*}
\int_{x_{i-1}}^{x_{i}}\int_{x_{0}}^{x_{i-1}}P(y)(t)\dt\dx & =\int_{x_{i-1}}^{x_{i}}\int_{x_{0}}^{x_{i-1}}y'(t)\dt\dx=\int_{x_{i-1}}^{x_{i}}\left(y(x_{i-1})-y(x_{0})\right)\dx\\
 & =\left(y(x_{i-1})-y_{0}\right)(x_{i}-x_{i-1}).
\end{align*}
Hence 
\begin{align*}
\int_{x_{0}}^{x_{n}}y(x)\dx-y_{0}(x_{n}-x_{0}) & =\sum_{i=1}^{n}\left(y(x_{i-1})-y_{0}\right)(x_{i}-x_{i-1})\\
 & \phantom{=}-(x_{n}-x_{0})\int_{x_{0}}^{x_{n}}y(x)\dx\sum_{i\ni a_{i}<0}a_{i}^{2}d_{i}+\frac{1}{2}(x_{n}-x_{0})^{2}\sum_{i=1}^{n}a_{i}^{2}f_{i},
\end{align*}
demonstrating that equation \eqref{eq:linear_system_part2} is satisfied
by \eqref{eq:linear_system_solution}.
\end{proof}
\begin{thm}
\label{thm:Fractal_transform} If $(y_{1},y_{2},\ldots,y_{n},A)$
is a solution of system \eqref{eq:linear_system_part1}, \eqref{eq:linear_system_part2}
and $f$ is an integrable function with $\int_{x_{0}}^{x_{n}}f(x)\dx=A$,
then
\begin{align*}
\mathcal{F}_{y_{0}}(f)(x_{k}) & =y_{k}\quad k=1,2,\ldots,n\\
\int_{x_{0}}^{x_{n}}\mathcal{F}_{y_{0}}(f)(x)\dx & =A
\end{align*}
\end{thm}

\begin{proof}
By proposition \ref{prop:basics2}, 
\[
\int_{x_{0}}^{x_{k}}P(f)(x)\dx=\sum_{i=1}^{k}\left[|a_{i}|d_{i}\int_{x_{0}}^{x_{n}}f(x)\dx+f_{i}(x_{i}-x_{i-1})\right].
\]
By definition, $\int_{x_{0}}^{x}P(f)(t)\dt=\mathcal{F}_{y_{0}}(f)(x)-y_{0}$,
so we have
\[
\mathcal{F}_{y_{0}}(f)(x_{k})=y_{0}+\sum_{i=1}^{k}\left[|a_{i}|d_{i}A+f_{i}(x_{i}-x_{i-1})\right]=y_{k}
\]
for $k=1,2,\ldots,n$. By proposition \ref{prop:basics3} and the
fact that $\sum_{i=1}^{n}\frac{a_{i}^{3}d_{i}}{|a_{i}|}=0$, 
\begin{align*}
\int_{x_{0}}^{x_{n}}\int_{x_{0}}^{x}P(f)(t)\dt\dx & =\sum_{i=1}^{n}\int_{x_{i-1}}^{x_{i}}\int_{x_{0}}^{x_{i-1}}P(f)(t)\dt\dx\\
 & \phantom{=}-(x_{n}-x_{0})\int_{x_{0}}^{x_{n}}f(x)\dx\sum_{i\ni a_{i}<0}a_{i}^{2}d_{i}+\frac{1}{2}(x_{n}-x_{0})^{2}\sum_{i=1}^{n}a_{i}^{2}f_{i}.
\end{align*}
Again using the fact that $\int_{x_{0}}^{x}P(f)(t)\dt=\mathcal{F}_{y_{0}}(f)(x)-y_{0}$,
\begin{align*}
\int_{x_{0}}^{x_{n}}\int_{x_{0}}^{x}P(f)(t)\dt\dx & =\int_{x_{0}}^{x_{n}}\left(\mathcal{F}_{y_{0}}(f)(x)-y_{0}\right)\dx\\
 & =\int_{x_{0}}^{x_{n}}\mathcal{F}_{y_{0}}(f)(x)\dx-y_{0}(x_{n}-x_{0})
\end{align*}
and 
\begin{align*}
\int_{x_{i-1}}^{x_{i}}\int_{x_{0}}^{x_{i-1}}P(y)(t)\dt\dx & =\int_{x_{i-1}}^{x_{i}}\left(\mathcal{F}_{y_{0}}(f)(x_{i-1})-y_{0}\right)\dx\\
 & =\left(\mathcal{F}_{y_{0}}(f)(x_{i-1})-y_{0}\right)(x_{i}-x_{i-1}).
\end{align*}
Hence 
\begin{align*}
\int_{x_{0}}^{x_{n}}\mathcal{F}_{y_{0}}(f)(x)\dx-y_{0}(x_{n}-x_{0}) & =\sum_{i=1}^{n}\left(\mathcal{F}_{y_{0}}(f)(x_{i-1})-y_{0}\right)(x_{i}-x_{i-1})\\
 & \phantom{=}-(x_{n}-x_{0})\int_{x_{0}}^{x_{n}}f(x)\dx\sum_{i\ni a_{i}<0}a_{i}^{2}d_{i}+\frac{1}{2}(x_{n}-x_{0})^{2}\sum_{i=1}^{n}a_{i}^{2}f_{i}.
\end{align*}
But $\mathcal{F}_{y_{0}}(f)(x_{i-1})=y_{i-1}$ so 
\begin{align*}
\int_{x_{0}}^{x_{n}}\mathcal{F}_{y_{0}}(f)(x)\dx & =y_{0}(x_{n}-x_{0})+\sum_{i=1}^{n}\left(y_{i-1}-y_{0}\right)(x_{i}-x_{i-1})\\
 & \phantom{=}-(x_{n}-x_{0})A\sum_{i\ni a_{i}<0}a_{i}^{2}d_{i}+\frac{1}{2}(x_{n}-x_{0})^{2}\sum_{i=1}^{n}a_{i}^{2}f_{i}\\
 & =A.
\end{align*}
\end{proof}
Theorems \ref{thm:FDE_soln_is_lin_sys_soln} and \ref{thm:Fractal_transform}
will form the cornerstone of our existence and uniqueness theorem
in the next section.

\section{SSDE's versus ODE's}

Ordinary differential equations of the form $y'=g(t,y)$ with initial
condition $y(t_{0})=y_{0}$ enjoy simple conditions for existence
and uniqueness of solutions (continuity of $g$ in $t$ and Lipschitz
continuity of $g$ in $y$), but, barring a known solution of the
ODE, do not admit simple calculation of $y(t)$ for any value other
than $t=t_{0}$. The study of numerical solutions of ODE's is dedicated
to the task of approximating such values \cite{brin,burden}. On the
other hand, SSDE's do not enjoy simple conditions for existence and
uniqueness of solutions (as best the authors can tell) but do admit
simple calculation of $y(t)$ for values other than $t=t_{0}$ under
certain conditions without having a solution of the SSDE in hand.
Despite the differences in form and calculability, proof of existence
and uniqueness of SSDE's proceeds along much the same lines as that
for ODE's, as seen in theorem \ref{thm:existence_uniqeness_1}.

The consequence of theorem \ref{thm:FDE_soln_is_lin_sys_soln} is
that, in certain instances, much can be determined about the solution(s)
of an SSDE, should any exist, before finding any solution. Under these
conditions, if system \eqref{eq:linear_system_part1}, \eqref{eq:linear_system_part2}
is inconsistent, then SSDE \eqref{eq:first_order} has no solution.
If system \eqref{eq:linear_system_part1}, \eqref{eq:linear_system_part2}
is consistent, then only certain sets of values $\left\{ \int_{x_{0}}^{x_{n}}y(x)\dx\right\} \cup\{y(x_{k}):k=1,2,\ldots,n\}$
are possible for solutions $y$. Equation \eqref{eq:transition} with
$y(0)=0$ is an example of an SSDE with infinitely many solutions,
one for each solution of system \eqref{eq:linear_system_part1}, \eqref{eq:linear_system_part2}.
As an example of an SSDE where the associated system \eqref{eq:linear_system_part1},
\eqref{eq:linear_system_part2} has exactly one solution, consider
the piecemealing of the form \eqref{eq:linear_piecemealing_FDE} on
$[1,4]$ given by the values in the following chart.
\begin{center}
\renewcommand{\arraystretch}{1.5}%
\begin{tabular}{|c|c|c|c|c|c|}
\multicolumn{1}{c}{$i$} & \multicolumn{1}{c}{$a_{i}$} & \multicolumn{1}{c}{$c_{i}$} & \multicolumn{1}{c}{$d_{i}$} & \multicolumn{1}{c}{$e_{i}$} & \multicolumn{1}{c}{$f_{i}$}\tabularnewline
\hline 
1 & $\frac{1}{3}$ & $0$ & $\frac{2}{3}$ & $\frac{2}{3}$ & $-\frac{22}{15}$\tabularnewline
\hline 
2 & $\frac{2}{3}$ & $0$ & $-\frac{1}{6}$ & $\frac{4}{3}$ & $\frac{17}{6}$\tabularnewline
\hline 
\end{tabular}
\par\end{center}

\noindent If \eqref{eq:first_order} with $y(1)=1$ has a solution
$y$, it is straightforward to compute $y(2)=\frac{23}{15}$, $y(4)=\frac{31}{5}$
and $\int_{1}^{4}y(x)\dx=9$ by solving system \eqref{eq:linear_system_part1},
\eqref{eq:linear_system_part2}.
\begin{thm}
\label{thm:existence-uniqueness}Let $P$ be a piecemealing of the
form \eqref{eq:linear_piecemealing_FDE} such that $\sum_{i=1}^{n}\frac{a_{i}^{3}d_{i}}{|a_{i}|}=0$.
If $\max\left\{ \frac{1}{2}|a_{i}d_{i}|:i=1,2,\ldots,n\right\} <1\}$
then there is a one-to-one correspondence between solutions of \eqref{eq:first_order}
on $[x_{0},x_{n}]$ and solutions of system \eqref{eq:linear_system_part1},
\eqref{eq:linear_system_part2}.
\end{thm}

\begin{proof}
Theorem \ref{thm:FDE_soln_is_lin_sys_soln} provides a function from
the set of solutions of SSDE \eqref{eq:first_order} to solutions
of system \eqref{eq:linear_system_part1}, \eqref{eq:linear_system_part2}.
Now suppose we have a solution $\{y_{1},y_{2},\ldots,y_{n},A\}$ of
system \eqref{eq:linear_system_part1}, \eqref{eq:linear_system_part2}.
Note that \eqref{eq:linear_system_part1} gives $y_{k}$ explicitly
in terms of $A$ for $k=1,2,\ldots,n$. Therefore this is the only
solution for this particular value of $A$. Now let $m=\max\left\{ \frac{1}{2}|a_{i}d_{i}|:i=1,2,\ldots,n\right\} <1\}$,
and let $g$ and $h$ be continuous functions such that ${\displaystyle \int_{x_{0}}^{x_{n}}g(x)\dx=\int_{x_{0}}^{x_{n}}h(x)\dx=A}$.
We wish to compare ${\displaystyle \max_{x\in[0,1]}\left|\mathcal{F}_{y_{0}}(g)(x)-\mathcal{F}_{y_{0}}(h)(x)\right|}$
with ${\displaystyle \max_{x\in[0,1]}\left|g(x)-h(x)\right|}$. First,
proposition \ref{prop:basics2} gives us that
\begin{align*}
\mathcal{F}_{y_{0}}(g)(x_{k})-\mathcal{F}_{y_{0}}(h)(x_{k}) & =\int_{x_{0}}^{x_{k}}\left(P(g)(t)-P(h)(t)\right)\dt\\
 & =\sum_{i=1}^{k}|a_{i}|d_{i}\int_{x_{0}}^{x_{n}}(g-h)(x)\dx\\
 & =0
\end{align*}
for $k=1,2,\ldots,n$. Hence 
\begin{align*}
\mathcal{F}_{y_{0}}(g)(x)-\mathcal{F}_{y_{0}}(h)(x) & =\int_{x_{0}}^{x}\left(P(g)(t)-P(h)(t)\right)\dt=\int_{x_{k}}^{x}\left(P(g)(t)-P(h)(t)\right)\dt
\end{align*}
 for any $k=1,\ldots,n$. In particular, given $1\leq i\leq n$, 
\begin{align*}
\mathcal{F}_{y_{0}}(g)(x)-\mathcal{F}_{y_{0}}(h)(x) & =\int_{x_{i-1}}^{x}d_{i}(g-h)\left(\frac{t-e_{i}}{a_{i}}\right)\dt\\
 & =a_{i}d_{i}\int_{\frac{x_{i-1}-e_{i}}{a_{i}}}^{\frac{x-e_{i}}{a_{i}}}(g-h)(u)\du
\end{align*}
 for $x_{i-1}\leq x\leq x_{i}$. Therefore
\begin{align*}
\max_{x\in[x_{i-1},x_{i}]}\left|\mathcal{F}_{y_{0}}(g)(x)-\mathcal{F}_{y_{0}}(h)(x)\right| & =\max_{x\in[x_{i-1},x_{i}]}\left|a_{i}d_{i}\int_{\frac{x_{i-1}-e_{i}}{a_{i}}}^{\frac{x-e_{i}}{a_{i}}}(g-h)(u)\du\right|\\
 & =\left|a_{i}d_{i}\right|\max_{x\in[x_{0},x_{n}]}\left|\int_{x_{0}}^{x}(g-h)(u)\du\right|,
\end{align*}
which implies
\[
\max_{x\in[x_{0},x_{n}]}\left|\mathcal{F}_{y_{0}}(g)(x)-\mathcal{F}_{y_{0}}(h)(x)\right|=2m\max_{x\in[x_{0},x_{n}]}\left|\int_{x_{0}}^{x}(g-h)(u)\du\right|.
\]
But $\int_{x_{0}}^{x_{n}}(g-h)(u)\du=0$, so by lemma \ref{lem:max}
\[
\max_{x\in[x_{0},x_{n}]}\left|\mathcal{F}_{y_{0}}(g)(x)-\mathcal{F}_{y_{0}}(h)(x)\right|\leq m\max_{x\in[x_{0},x_{n}]}\left|g(x)-h(x)\right|.
\]
Hence, for this value $A$, $\mathcal{F}_{y_{0}}$ is contractive
on 
\[
Z_{A}=\left\{ f:[0,1]\to\mathbb{R}|f\text{ is continuous and }\int_{x_{0}}^{x_{n}}f(x)\dx=A\right\} 
\]
 under the max norm. By theorem \ref{thm:Fractal_transform}, $Z_{A}$
is closed under $\mathcal{F}_{y_{0}}$. Since $Z_{A}$ is complete
with respect to uniform convergence, the contraction mapping principle
guarantees a unique fixed point of $\mathcal{F}_{y_{0}}$ in $Z_{A}$.
This fixed point, then, is the unique solution of the equation $y=\mathcal{F}_{y_{0}}(y)$,
the integral equation equivalent of \eqref{eq:first_order}, within
$Z_{A}$.
\end{proof}

\section{Higher Order SSDE's}

Higher order SSDE's can be defined analogous to \eqref{eq:first_order}
in the following way. A differential equation of the form 
\begin{align}
y^{(n)} & =P(y)\label{eq:nth_order}\\
y(x_{0}) & =y_{0},\ y'(x_{0})=y'_{0},\ldots,y^{(n-1)}(x_{0})=y_{0}^{(n-1)}\nonumber 
\end{align}
where $P$ is a piecemealing on $I=[x_{0},x_{n}]$ is called an $n^{th}$
order SSDE. In fact, the genesis of this research lies in the question
of existence, uniqueness, and computability of the solution of a second
order SSDE related to cam design. In this original formulation, we
seek a transition function that reaches its maximum acceleration quickly,
holds that acceleration for a time, then transitions quickly to its
minimum acceleration, and holds that for a time before returning to
acceleration zero. The transitions between states of constant acceleration
are to be linear transformations of the entire transition function.
Because acceleration and force are proportional, optimizing for a
low maximum acceleration amounts to optimizing for a cam subject to
low forces. Further, because transition functions are smooth, the
third derivative, also called the jerk, (and all higher derivatives)
will be bounded, unlike many piecewise obtained cam profiles. More
specifically, and mathematically, we seek a solution of the second
order SSDE
\begin{align*}
y'' & =P(y)\\
y(0) & =0,\ y'(0)=0
\end{align*}
 where $P(y)$ is a piecemealing of the form \eqref{eq:linear_piecemealing_FDE}
on $[0,1]$ given by the values in the following chart.
\begin{center}
\renewcommand{\arraystretch}{1.5}%
\begin{tabular}{|c|c|c|c|c|c|}
\multicolumn{1}{c}{$i$} & \multicolumn{1}{c}{$a_{i}$} & \multicolumn{1}{c}{$c_{i}$} & \multicolumn{1}{c}{$d_{i}$} & \multicolumn{1}{c}{$e_{i}$} & \multicolumn{1}{c}{$f_{i}$}\tabularnewline
\hline 
1 & $\frac{1}{6}$ & $0$ & $6$ & $0$ & $0$\tabularnewline
\hline 
2 & $\frac{1}{6}$ & $0$ & $0$ & $\frac{1}{6}$ & $6$\tabularnewline
\hline 
3 & $\frac{1}{6}$ & $0$ & $-6$ & $\frac{1}{3}$ & $6$\tabularnewline
\hline 
4 & $\frac{1}{6}$ & $0$ & $-6$ & $\frac{1}{2}$ & $0$\tabularnewline
\hline 
5 & $\frac{1}{6}$ & $0$ & $0$ & $\frac{2}{3}$ & $-6$\tabularnewline
\hline 
6 & $\frac{1}{6}$ & $0$ & $6$ & $\frac{5}{6}$ & $-6$\tabularnewline
\hline 
\end{tabular}
\par\end{center}

\begin{figure}
\begin{centering}
(a) \raisebox{-.75in}{\includegraphics[width=1.75in]{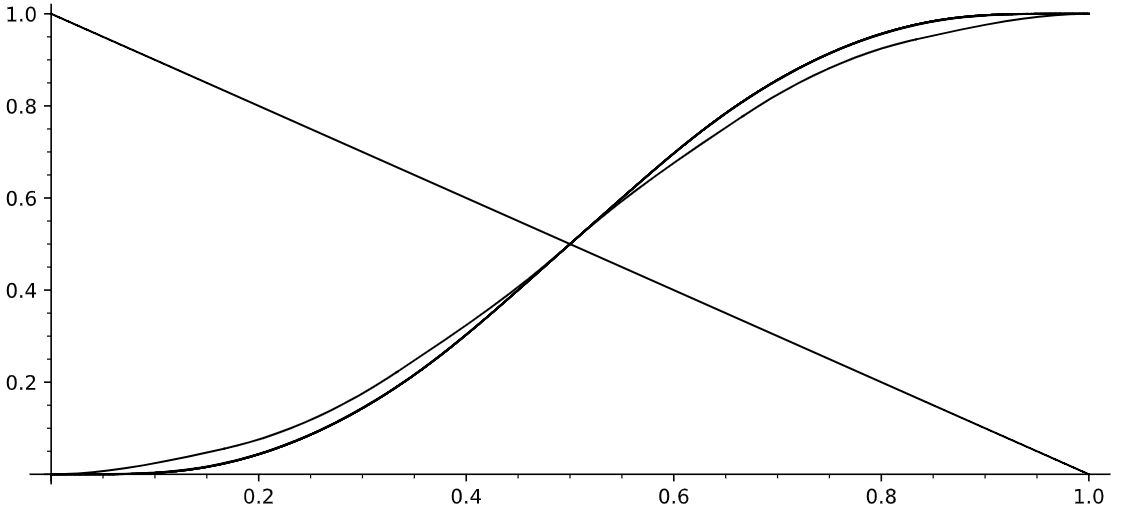}}
(b) \raisebox{-.75in}{\includegraphics[width=1.75in]{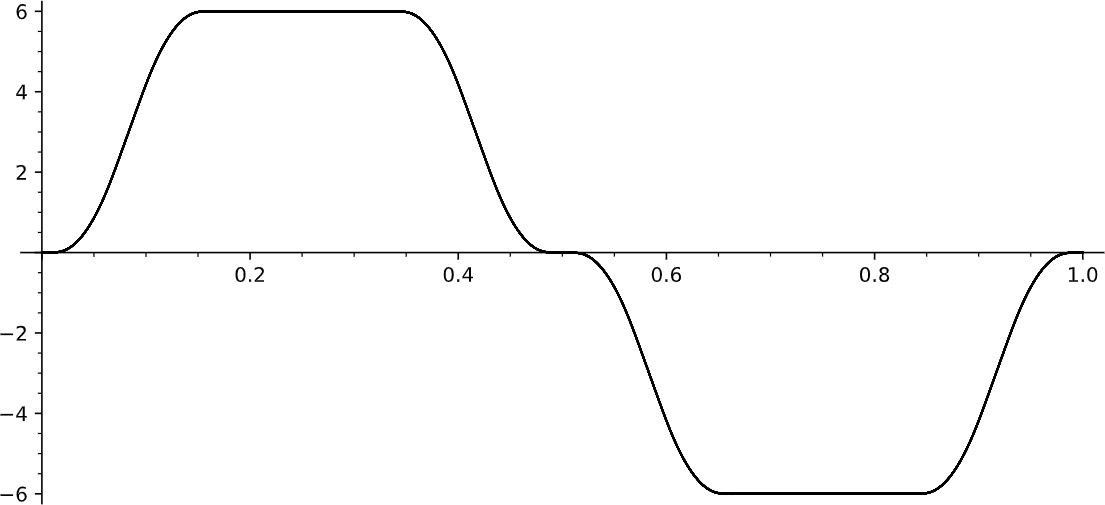}
}
\par\end{centering}
\caption{(a) iterates of $S$ beginning with $f_{0}(x)=1-x$, (b) $P(f_{4})$
and $f_{4}''$} \label{fig:second_order_FDE}
\end{figure}

\noindent We state evidence without proof that this SSDE has a unique
solution and that the solution is a transition function. Setting $S(g)(x)=\int_{0}^{x}\int_{0}^{t}g(u)\du\dt$
and $f_{0}(x)=1-x$, we compute $f_{1}(x)=S(f_{0})(x)$, $f_{2}(x)=S(f_{1})(x)$,
and so on through $f_{4}(x)$. Results of this calculation appear
in Figure \ref{fig:second_order_FDE}. Note that the graphs of $P(f_{4})$
and $f_{4}''$ are indistinguishable. 

\end{document}